\documentclass[11pt]{article}
\usepackage[colorlinks=true,
linkcolor=webgreen,
filecolor=webbrown,
citecolor=webgreen]{hyperref}

\usepackage{amssymb,psfig,epsfig}

\definecolor{webgreen}{rgb}{0,.5,0}
\definecolor{webbrown}{rgb}{.6,0,0}

\newtheorem{theorem}{Theorem}

\def\slfrac#1#2{\hbox{\kern.1em %
 \raise.5ex\hbox{\the\scriptfont0 #1}\kern-.11em %
 /\kern-.15em\lower.25ex\hbox{\the\scriptfont0 #2}}}

\newcommand{\eqn}[1]{(\ref{#1})}
\newcommand{\hsp}{\hspace*{\parindent}}

\newcommand{\eeq}{\end{equation}}
\newcommand{\beql}[1]{\begin{equation}\label{#1}}
\newcommand{\bsq}{{\vrule height .9ex width .8ex depth -.1ex }}

\newcommand{\RR}{\mathbb R}
\newcommand{\ZZ}{\mathbb Z}
\newcommand{\NN}{\mathbb N}
\newcommand{\QQ}{\mathbb Q}
\newcommand{\af}{\alpha}

\makeatletter
\def\@sect#1#2#3#4#5#6[#7]#8{\ifnum #2>\c@secnumdepth
     \def\@svsec{}\else
     \refstepcounter{#1}\edef\@svsec{\csname the#1\endcsname.\hskip .75em }\fi
     \@tempskipa #5\relax
      \ifdim \@tempskipa>\z@
        \begingroup #6\relax
          \@hangfrom{\hskip #3\relax\@svsec}{\interlinepenalty \@M #8\par}%
        \endgroup
       \csname #1mark\endcsname{#7}\addcontentsline
         {toc}{#1}{\ifnum #2>\c@secnumdepth \else
                      \protect\numberline{\csname the#1\endcsname}\fi
                    #7}\else
        \def\@svsechd{#6\hskip #3\@svsec #8\csname #1mark\endcsname
                      {#7}\addcontentsline
                           {toc}{#1}{\ifnum #2>\c@secnumdepth \else
                             \protect\numberline{\csname the#1\endcsname}\fi
                       #7}}\fi
     \@xsect{#5}}
\def\@begintheorem#1#2{\it \trivlist \item[\hskip \labelsep{\bf #1\ #2.}]}

\def\plain{plain}\ifx\fmtname\plain\csname fi\endcsname
     
     \input docstrip
     \preamble

     Do not distribute the stripped version of this file.
     The checksum in the header refers to the documented version.

     \endpreamble
     \generateFile{here.sty}{t}{\from{here.doc}{}}
     \endinput
\fi
\ifcat a\noexpand @\let\next\relax\else\def\next{%
    \documentstyle[here,doc]{article}\MakePercentIgnore}\fi\next
\ifx\@Hxfloat\@Hundef\else\expandafter\endinput\fi
\let\@Hxfloat\@xfloat
\def\@xfloat#1[{\@ifnextchar{H}{\@HHfloat{#1}[}{\@Hxfloat{#1}[}}
\def\@HHfloat#1[H]{%
\expandafter\let\csname end#1\endcsname\end@Hfloat
\vskip\intextsep\vbox\bgroup\def\@captype{#1}\parindent\z@
\ignorespaces}
\def\end@Hfloat{\egroup\vskip \intextsep}

\makeatother

\begin{document}
\begin{center}
{\Large {\bf On the Existence of Similar Sublattices}} \\
\vspace{1.5\baselineskip}
{\em J. H. Conway} \\
Mathematics Department \\
Princeton University \\
Princeton, NJ 08540 \\
\vspace{1\baselineskip}
{\em E. M. Rains} and {\em N. J. A. Sloane} \\
Information Sciences Research \\
AT\&T Shannon Lab \\
Florham Park, NJ 07932-0971 \\
\vspace{1.5\baselineskip}
\end{center}
November 17, 1998; revised October 8, 1999.
A slightly different version of this paper appeared in the
{\em Canadian Jnl. Math.}, {\bf 51} (1999), pp. 1300--1306.
Pdf version prepared May 7, 2000.
\begin{center}
\vspace{2.0\baselineskip}
DEDICATED TO H. S. M. COXETER \\
\vspace{2.0\baselineskip}
{\bf ABSTRACT}
\vspace{.5\baselineskip}
\end{center}
\setlength{\baselineskip}{1.5\baselineskip}

Partial answers are given to two questions.
When does a lattice $\Lambda$ contain a sublattice $\Lambda '$
of index $N$ that is geometrically similar to $\Lambda$?
When is the sublattice ``clean'', in the sense that the boundary of the
% $\Lambda '$ 
Voronoi cells for $\Lambda '$ do not intersect $\Lambda$?

\clearpage
\thispagestyle{empty}
\setcounter{page}{1}

\section{Introduction}
\hsp
A {\em similarity} $\sigma$ of norm $c$ is a linear map from
$\RR^n$ to $\RR^n$ such that
$\sigma u \cdot \sigma v = c ~ u \cdot v$ for $u,v \in \RR^n$.
Let $\Lambda$ be an $n$-dimensional rational lattice, i.e.
$u \cdot v \in \QQ$ for $u,v \in \Lambda$.
A sublattice $\Lambda ' \subseteq \Lambda$ is {\em similar}
to $\Lambda$ if $\sigma ( \Lambda ) = \Lambda '$ for some similarity
$\sigma$ of norm $c$.
We also call $\sigma$ a {\em multiplier of norm $c$} for $\Lambda$.
The index $N = [\Lambda : \Lambda' ]$ is $c^{n/2}$, so if $n$ is odd
$c$ must be a square, say $c =a^2$, and we could take $\sigma$ to be scalar
multiplication by $a$.
In other words the norms of similarities of odd-dimensional
lattices are precisely the integral squares.
Henceforth we will assume that $n=2k$ is even.

Multipliers of small norm, especially 2 (also called
``norm-doubling maps'') are useful for recursive constructions
of lattices (\cite{SPLAG}, Chap. 8).
If the root lattice $E_6$ had a norm-doubling map $\sigma$, then the
``$u$, $u+v$'' construction\footnote{Take the lattice consisting of the vectors
$(u,u+v)$ for $u \in E_6$, $v \in \sigma (E_6)$.}
would produce a denser 12-dimensional lattice than the Coxeter-Todd lattice $K_{12}$.
However, some years ago W. M. Kantor and N.J.A.S. showed by
direct search that no such map exists.
This result now follows from Theorem~\ref{th2}.

The question of the existence of multipliers of given norm
arose recently in constructing ``multiple description''
vector quantizers \cite{MD,VSS}.
In an ordinary vector quantizer an $n$-dimensional lattice $\Lambda$ is specified,
and successive $n$-tuples $(x_1, \ldots, x_n) \in \RR^n$ are
replaced by the closest lattice points (cf. \cite[Chapter 2]{SPLAG}).
In a multiple description scheme we also choose a number $N$ and a labeling
$$u \in \Lambda \mapsto (l(u), r(u)) \in \ZZ \times \ZZ$$
such that $| l^{-1} (i) \cap r^{-1} (j) | \le 1$, $| l^{-1} (i) | \le N$,
$|r^{-1} (j) | \le N$ for all $i,j \in \ZZ$.
The numbers $l(u)$ and $r(u)$ are transmitted over different
channels.
If both numbers are received then $u$ is uniquely determined, but if only one number is received (and the other lost) then $u$ is determined to within a small region of $\RR^n$
(and the goal is to choose $l$ and $r$ so that this region is as small as possible).
The method proposed in \cite{MD,VSS} for constructing such labelings makes use of a
sublattice $\Lambda'$ that has index $N$ in $\Lambda$ and is similar to $\Lambda$.
For this application it is also of interest to know when the boundary of the Voronoi cell of the
sublattice $\Lambda'$ does not contain any points of $\Lambda$:
we call such sublattices ``clean''.

In Section 2 we give several results about the existence of similar
sublattices,
and in Section 3 we give a partial answer to the existence of clean
sublattices in the two-dimensional case.

%It is worth mentioning that in some cases
%Siegel's ``mass formula'' could be used to
%determine how many similar sublattices exist of given norm.

The only references we have found which treat the
first problem are Baake and Moody \cite{BM1}, \cite{BM2},
which are concerned with lattices (and more
general $\ZZ$-modules related to quasicrystals)
in dimensions 1 to 4.
These authors use techniques from ideal theory
and quaternion algebras to
enumerate similar substructures of given index.

A related problem has been studied in the crystallographic literature
\cite{B1}, \cite{B2}, \cite{B3}, \cite{B4}:
given a lattice $\Lambda$ (or more generally a $\ZZ$-module) in $\RR^n$,
when does there exist an isometry $\sigma$ such that the
``coincidence site sublattice'' $\Lambda ' = \Lambda \cap \sigma (\Lambda )$ has finite
index in $\Lambda$?
This is a somewhat different problem, since $\Lambda '$ need not be
similar to $\Lambda$, nor can every similar sublattice of $\Lambda$ be obtained
in this way.

We discovered the above references by accident.
Using a computer we found that 
the lattice $A_4$ has multipliers of norm $c$
precisely when $c$ is one of the numbers
$$1, 4, 5, 9, 11, 16, 19, 20, 25, 29, 31, 36, \ldots$$
The same sequence appears in \cite{B1}\footnote{Found with the help of \cite{EIS}, where it is Sequence 
\htmladdnormallink{A31363}{http://www.research.att.com/cgi-bin/access.cgi/as/njas/sequences/eisA.cgi?Anum=031363}.}, as the indices of coincidence
site sublattices in a certain three-dimensional quasicrystal.
\cite{B1} identifies these numbers as those
positive integers in which all primes congruent to 2 or 3 ($\bmod~5$) appear to
an even power.
As Theorem~\ref{th2} shows, this is the same as our sequence.
Although this can hardly be a coincidence, we do not at present
see a direct connection between the $A_4$ and quasicrystal problems.

Two papers by Chapman \cite{Chap1}, \cite{Chap2} consider a different,
though again related, problem concerning sublattices of $\ZZ^n$.

\section{The existence of similar sublattices}
\hsp
Let $\Lambda$ be a rational $2k$-dimensional lattice with Gram matrix $A$,
and let $c \in \NN$.  We wish to know if  $\Lambda$ has a sublattice
$\Lambda '$ such that $\sigma ( \Lambda ) = \Lambda '$
for some similarity $\sigma$ of norm $c$.
The existence of $\Lambda '$ can be determined (in principal)
by searching through  $\Lambda$ to see if it contains a set of
vectors with Gram matrix $cA$.  For small $k$ and $c$ this
is quite feasible.  We know of no other
method that will always succeed.

By using the rational invariants of  $\Lambda$ we can obtain
a necessary condition for $\Lambda '$ to exist, which in some
cases is also sufficient.
The Hilbert symbol (\cite{Kit}, \cite{Wat})
provides a convenient way to specify this condition.

For a rational number $r > 0$ we write $(r)$ for the fractional
ideal $r\QQ$.  A lattice $\Lambda$ is {\em $(r)$-maximal}
if $\Lambda$ is maximal with respect to the property that
$u \cdot u \in (r)$ for all $u \in \Lambda$ (\cite{Kit}, \cite{OM}).
The importance of this concept stems from the result (\cite{OM}, Section 102:3)
that the $(r)$-maximal lattices in a rational class
form a single genus.
We also say that
$\Lambda$ is {\em unigeneric} if it is unique in its genus.

\begin{theorem}\label{th1}
A necessary condition for a $2k$-dimensional lattice $\Lambda$ to have a
multiplier of norm $c$ is that the Hilbert symbol
\beql{Eq1}
(c,(-1)^k \det \Lambda)_p =1
\eeq
for all primes $p$ dividing $2c \det \Lambda$.
If $\Lambda$ is unigeneric and $(r)$-maximal for some $r \in \QQ$
then this condition
is also sufficient.
\end{theorem}

\paragraph{Proof.}
If $\sigma$ is a multiplier of norm $c$ for $\Lambda$, then
$\Lambda ' = \sigma (\Lambda )$ and $\Lambda$ are rationally equivalent,
hence equivalent over the $p$-adic rationals for all $p$, and so have the same
Hasse-Minkowski invariant $\epsilon_p$ for all $p$.
The $p$-adic Hasse-Minkowski invariants for $\Lambda$ and $\Lambda '$ differ by a factor of $(c, (-1)^k \det \Lambda )_p$
(\cite{Wat}, p. 46; \cite{Kit}, Theorem 3.4.2).
Since this invariant is 1 if $p$ does not divide $2c \det \Lambda$,
the first assertion
follows.
Conversely, if \eqn{Eq1} holds for all $p$ then $\Lambda$ and the rescaled
lattice ${\sqrt{c}} \Lambda$ are rationally equivalent.
For $u \in {\sqrt{c}} \Lambda$, $u \cdot u \in (cr) \subseteq (r)$,
so ${\sqrt{c}} \Lambda$ is contained in some maximal $(r)$-lattice $M$, say.
By \cite{OM}, Section 102:3, $M$ and $\Lambda$ are in the same
genus, and since $\Lambda$ is unigeneric, $M$ is in
the same class as $\Lambda$.
Hence $\Lambda$ has a sublattice in the same class as ${\sqrt{c}} \Lambda$.  ~~~$\bsq$

Note that it is not enough for $\Lambda$ to be unigeneric for
the condition of the theorem to be sufficient.  The lattice
with Gram matrix
$\left[ {1 \atop 0} ~ {0 \atop 4} \right]$
is unigeneric and rationally equivalent to
$\left[ {2 \atop 0} ~ {0 \atop 8} \right]$,
but does not have a similarity of norm 2.

Many familiar lattices of small determinant and dimension are unigeneric (see \cite{LDLI}),
and are often $(1)$-maximal (if they contain
vectors of odd norm) or $(2)$-maximal
(if they contain only vectors of even norm).
In these cases Theorem \ref{th1} provides the answer
to our first question, as the following theorem illustrates.
The straightforward proof is omitted.

\begin{theorem}\label{th2}
The lattices $\ZZ^2$, $A_2$, $A_4$, $\ZZ^6$, $E_6$ have multipliers of norm $c$ just for the following values:
$$
\begin{array}{l}
\mbox{$\ZZ^2$ or $\ZZ^6$: $c = r^2 + s^2$, $r, s \in \ZZ$ $($Sequence
\htmladdnormallink{A001481}{http://www.research.att.com/cgi-bin/access.cgi/as/njas/sequences/eisA.cgi?Anum=001481}$)$; i.e.} \\ [+.1in]
\mbox{~~~~primes $\equiv 3 ~( \bmod ~4)$ appear to even powers in $c$,} \\ [+.2in]
\mbox{$A_2$ or $E_6$: $c = r^2 - rs + s^2$, $r, s \in \ZZ$ $($Sequence
\htmladdnormallink{A003136}{http://www.research.att.com/cgi-bin/access.cgi/as/njas/sequences/eisA.cgi?Anum=003136}$)$; i.e.} \\ [+.1in]
\mbox{~~~~primes $\equiv 2 ~( \bmod ~3)$ appear to even powers in $c$,} \\ [+.2in]
\mbox{$A_4$: $c = r^2 + rs -s^2$, $r, s \in \ZZ$ $($Sequence
\htmladdnormallink{A031363}{http://www.research.att.com/cgi-bin/access.cgi/as/njas/sequences/eisA.cgi?Anum=031363}$)$; i.e.} \\ [+.1in]
\mbox{~~~~primes $\equiv \pm 2 ~( \bmod~5)$ appear to even powers in $c$.}
\end{array}
$$
\end{theorem}

In some cases explicit similarities are easily found.
For $\ZZ^2$ and $A_2$ we use complex coordinates and take $\sigma$
to be multiplication by $r+si$ and $r+ \omega s$ respectively,
where $\omega  = e^{2 \pi i/3}$.
For $E_6$ we use three complex coordinates (\cite{SPLAG}, p. 126, Eq. (120))
and again multiply by $r+ s \omega$.

For $A_4$ some similarities can be found using the element
$$
\alpha = \left[
\begin{array}{ccccc}
0 & 1 & 0 & 0 & 0 \\
0 & 0 & 1 & 0 & 0 \\
0 & 0 & 0 & 1 & 0 \\
0 & 0 & 0 & 0 & 1 \\
1 & 0 & 0 & 0 & 0
\end{array}
\right]
$$
of $Aut (A_4)$.
It can be shown that $\sigma = a_1 \af + a_2 \af^2 + a_3 \af^3 + a_4 \af^4$ is a similarity for $A_4$ of norm
$$
\begin{array}{cccc}
\frac{1}{2} \, [a_1 & a_2 & a_3 & a_4] \\
&&&\\
&&&\\
&&&
\end{array}
\left[
\begin{array}{rrrr}
2 & -1 & 0 & 0 \\
-1 & 2 & -1 & 0 \\
0 & -1 & 2 & -1 \\
0 & 0 & -1 & 2
\end{array}
\right]~
\left[
\begin{array}{c}
a_1 \\ a_2 \\ a_3 \\ a_4
\end{array}
\right]
$$
provided
$$
\begin{array}{cccc}
[a_1 & a_2 & a_3 & a_4] \\
&&&\\
&&&\\
&&&
\end{array}
\left[
\begin{array}{rrrr}
0 & 1 & -1 & -1 \\
1 & 0 & 1 & -1 \\
-1 & 1 & 0 & 1 \\
-1 & -1 & 1 & 0
\end{array}
\right]~
\left[
\begin{array}{c}
a_1 \\ a_2 \\ a_3 \\ a_4
\end{array}
\right] =0 ~.
$$
This gives similarities of norms $1, 5, 11, \ldots$ but not $19, 29, \ldots$~.
We do not know a simple way to find the other similarities.
Of course analogous similarities can be found for any cyclotomic lattice.

\begin{theorem}\label{th3}
The lattices $\ZZ^{4m}$, $D_{4m}$ and $D_{4m}^+$
$(m \ge 1)$,
$E_8$, $K_{12}$, the Barnes-Wall lattice $BW_{16}$ and the Leech
lattice $\Lambda_{24}$ have multipliers of every norm.
\end{theorem}

\paragraph{Remark.}
Baake and Moody \cite{BM2} establish this for
$\ZZ^4$ and $D_4$ and also give a Dirichlet
generating function for the number
of similar sublattices of given index.

\paragraph{Proof.}
For $\ZZ^{4m}$, $D_{4m}$, $D_{4m}^+$ we represent
the vectors by $m$ Hurwitz integral quaternions;
then right multiplication by $q = r+si+tj+uk$ is a similarity of norm
$|q|^2 = r^2 +s^2 + t^2 +u^2$.

We write the vectors of $\Lambda_{24}$ in $4 \times 6$ MOG coordinates
\cite{SPLAG} and convert each of the six columns to a quaternion according to the scheme
$$
\begin{array}{|cc|cc|cc|} \hline
1 & 1 & 1 & 1 & 1 & 1 \\
k & j & k & j & k & j \\ \hline
i & k & i & k & i & k \\
j & i & j & i & j & i \\ \hline
\end{array} ~.
$$
In this form left or right multiplication by $i, j$ or $\omega = \frac{1}{2} (-1 + i+j+k)$ are all
automorphisms of $\Lambda_{24}$, as is the column
permutation $(1,2) (3,4) (5,6)$.
Then right multiplication by $q$ is a similarity of $\Lambda_{24}$ of norm $|q|^2$.

We define $BW_{16}$ to be the sublattice of $\Lambda_{24}$ in which the
last two columns of the MOG are zero, and again use $q$.
Finally, we define $K_{12}$ to be the sublattice of $\Lambda_{24}$ consisting of vectors
$$
\begin{array}{|cc|cc|cc|} \hline
a & b & c & d & e & f \\
A & B & C & D & E & F \\ \hline
A & B & C & D & E & F \\
A & B & C & D & E & F \\ \hline
\end{array} ~,
$$
which we associate with the three-dimensional quaternionic vector
\beql{Eq2}
(a+ bi + \sqrt{3} Aj + \sqrt{3} Bk , \ldots, e+ fi + \sqrt{3} Ej + \sqrt{3} Fk ) ~.
\eeq
Then right multiplication of \eqn{Eq2} by $r+ si + \sqrt{3} tj + \sqrt{3} uk$ defines a similarity
of norm $r^2 + s^2 + 3t^2 + 3u^2$.
Since the latter form represents $1, \ldots, 15$, by the ``15-Theorem''
of J.H.C. and W. A. Schneeberger (cf. \cite{SPLAG}) it
represents all numbers, and the proof is complete.
(Of course the universality of this form
was already known: it appears in \cite{Ram}.)~~~$\bsq$

Another easy consequence of Theorem \ref{th1} is:
\begin{theorem}\label{th4}
A necessary condition for $\Lambda$ to have a norm-doubling map is that
$\dim \Lambda$ be even and that all primes $\equiv \pm 3 ~( \bmod~8)$ appear
to even powers in $\det \Lambda$.
If $\Lambda$ is unigeneric
and $(r)$-maximal for some $r \in \QQ$
then this condition is also sufficient.
\end{theorem}

\section{The existence of clean sublattices in the two-dimensional case}
\hsp
The Voronoi cell of a two-dimensional lattice is either a hexagon or a rectangle (\cite{LDL6}, Fig.~1).
We assume that the lattice is generated by 1 and an imaginary quadratic integer.
Similar arguments could be applied to more general two-dimensional lattices but
the answers would be much more complicated.

We first consider a lattice $\Lambda$ with a hexagonal
Voronoi cell, generated say by 1 and $\omega = (-1 + \sqrt{-N}) /2$,
$N \equiv 3$ $(\bmod~4)$.
A similarity $\sigma$ of norm $c$ is represented by multiplication by $\af = a+ b \omega$,
$a,b \in \ZZ$, with $c = | \af |^2 = a^2 - ab + (N+1) b^2/4$.
We begin with the case $N=3$, the hexagonal lattice (a rescaled version of $A_2$).

\begin{theorem}\label{th5}
For the hexagonal lattice generated by $1$ and $\omega = e^{2 \pi i/3}$,
multiplication by $\af = a + b \omega$ yields a clean sublattice if and only if
$\af \theta$ is a primitive element\footnote{An element $r+ s \omega \in \ZZ [\omega]$ is primitive
if and only if $gcd (r,s) =1$.} of $\ZZ [\omega]$,
where $\theta = \omega - \bar{\omega} = \sqrt{-3}$.
There is a clean sublattice of index $c$
if and only if $c$ is a product of primes
$\equiv 1$ $(\bmod ~3)$.
\end{theorem}

\paragraph{Proof.}
The Voronoi cell of $\Lambda$ is a regular hexagon,
and all edges are equivalent, so it is enough to consider say the left-hand edge $L$.
This edge is the middle third of the line $M$ from $\omega$ to $\bar{\omega}$.
The lattice
$\Lambda ' = \af \Lambda$ is clean if and only if there is no point of $\Lambda$ on the line
$\af L$.
Since $\Lambda$ is a lattice, if there is a point in the interior of $M$ then there is a point on $L$.
So $\Lambda'$ is clean if and only if there is no point of $\Lambda$ in the interior of $\af M$, i.e. if and only if $\af \theta$ is primitive.
The second assertion follows easily from the fact that the numbers primitively represented
by $a^2 - ab + b^2$ are of the form $3^\epsilon$ times a product of distinct primes $\equiv 1$
$(\bmod~3)$, where $\epsilon =0$ or 1.~~~$\bsq$

We state the result for the general case without proof.
The argument is similar to the above, but one must consider all sides of the Voronoi cell.
\begin{theorem}\label{th6}
For the hexagonal-type lattice generated by $1$
and $\omega = (1+ \sqrt{-N})/2$, the similarity defined by multiplication by $\af = a+b \omega$ yields a clean sublattice
if and only if 
\begin{itemize}
\item[(i)]
$\af \theta$ is primitive, where $\theta = \omega - \bar{\omega} = \sqrt{-N}$,
\item[(ii)]
there is an odd number $k$ dividing $N+1$ such that $\af (N- \theta ) / (2k)$ is integral and primitive, and
\item[(iii)]
there is an odd number $k$ dividing $N+1$ such that $\af (N+ \theta )/(2k)$ is integral and primitive.
\end{itemize}
\end{theorem}

For the Kleinian lattice generated by 1 and $\omega = (1+ \sqrt{-7} )/2$, for example,
$k=1$ is the only possibility, and the theorem states that $\af \Lambda$ is
clean if and only if $\af \theta$, $\af \omega \theta$ and $\af \bar{\omega} \theta$ are all
primitive.
Equivalently, $(a+b \omega) \Lambda$ is primitive if and only if $a$ is odd,
$b$ is even and $gcd (a,b) =1$.

Similar arguments also give the result for lattices with rectangular Voronoi cells.
Again we omit the proof.

\begin{theorem}\label{th7}
For the rectangular-type lattice generated by $1$ and $\theta = \sqrt{-N}$, $N \ge 1$,
the similarity defined by multiplication by $\af = a+b \sqrt{-N}$ yields
a clean sublattice if and only if $\af \bar{\af} = a^2 + Nb^2$ (or equivalently
$a+Nb$) is odd.
\end{theorem}

\paragraph{Remarks.}
M. Baake has pointed out to us that by combining Theorems \ref{th5}
and \ref{th7} with the results of \cite{B1} we can say that
a sublattice of $A_2$ (or $\ZZ^2$) is clean if and only if it is a
coincidence site sublattice.

It would be nice to know what happens in higher dimensions.
What is the analogue of Theorem \ref{th5} for $D_4$ or $E_8$, for instance?

\paragraph{Acknowledgments.}
We thank H.-G. Quebbemann for suggesting that the similar sublattice problem
is best handled via the Hilbert symbol, and for pointing out an error
in the first version of this paper.
We also thank
V. A. Vaishampayan and S. D. Servetto for many conversations
concerning vector quantization,
and M. Baake and J. Martinet for comments on the manuscript.
R. Schulze-Pillot drew our attention to the Chapman references.

\end{document}